\documentclass[12pt]{article}

\newcommand{\roim}[1]{{\mathrm{R}}{#1}_*}

\newcommand{\shl}{{\cal L}}
\newcommand{\shf}{{\cal F}}

\newcommand{\BDCC}{{{\bf D}^b_{c}}}

\title{ 
Addendum to the paper \\
``Characteristic Cycles 
of Perverse Sheaves  \\
and Milnor Fibers"}

\date{}

\begin{document}  

\maketitle
\pagestyle{plain}

\bigskip

\begin{center}
{\large   Philibert Nang and Kiyoshi Takeuchi }
\end{center}

\begin{center}
{ Institute of Mathematics, 
University  of Tsukuba, 
1-1-1, Tennodai, Tsukuba, Ibaraki, 305-8571, JAPAN, 
\\ 
E-mail:  takechan@math.tsukuba.ac.jp
\\ 
}
\end{center}

\bigskip
\noindent We note that 
Theorem 5.4 of the paper ``Characteristic cycles 
of perverse sheaves  
and Milnor fibers" (by P. Nang and K. Takeuchi, 
published online in Math. Zeitschrift, 
October 15, 2004) holds for ``any" non-zero complex number $a$.  
Indeed, using the notations of this paper, 
let us change the definition of 
$\shf^{\cdot}\in \BDCC (X)$ in the proof of Theorem 5.4 as 
$$\shf^{\cdot}:=
\tau^{\leq n-1}\roim{j}\widetilde{\shl_{a}}.$$ 
\noindent Then $\shf^{\cdot}$ is a perverse sheaf on $X$ since 
$j_{0}:X\backslash V\hookrightarrow X\backslash 
\left\{ 0\right\} $ is a Stein map.  Moreover by Lemma 5.5 
we have $\chi_x (\shf^{\cdot})=0$ at any point 
$x\in V$.  Therefore, without assuming any technical 
condition on the complex number $a\not= 0$,  
the proof proceeds exactly in the 
same way as before until we get the result 
$$N_a \leq (-1)^n\{\chi (V\cap L\cap B_{\varepsilon})-1\},$$ 
\noindent where the right-hand-side 
is noting but the top-dimensional Betti number 
of the complex 
link $V\cap L\cap B_{\varepsilon}$ of $V$, 
i.e. the number of spheres appearing in the bouquet 
decomposition of the complex link.   In the same way 
we can also completely remove the assumption of the 
complex number $a\not= 0$ in 
Proposition 6.4.17 of A. Dimca's book 
``Sheaves in Topology"  Springer 2004 
(which slightly generalized 
our previous Theorem 5.4).

\end{document}